\documentclass[12pt,oneside,a4paper,reqno]{article}
\usepackage{latexsym,amsxtra,amsthm, amsmath, amsfonts, amscd, amssymb, verbatim}
\usepackage{amsfonts}
\usepackage[portrait, top=2.5 cm, bottom=3cm, left=3.5cm, right=2cm] {geometry} 

\usepackage{varioref}
\usepackage{fancyhdr}
\usepackage{longtable}
\usepackage{fancybox}
\usepackage{color, graphicx}
\usepackage{verbatim}
\usepackage{setspace}

\usepackage[latin1]{inputenc}
\usepackage[english]{babel}

\begin{document}

\begin{center}
\textbf{COUNTING RATIONAL POINTS ON SMOOTH CUBIC CURVES}
\end{center}

\begin{center}
\textbf{Manh Hung Tran}
\end{center}

\begin{center}
Department of Mathematical Sciences, Chalmers University of Technology.
\end{center}

\begin{center}
\textbf{ABSTRACT}
\end{center}
We use a global version of Heath-Brown's $p-$adic determinant method developed by Salberger to give upper bounds for the number of rational points of height at most $B$ on non-singular cubic curves defined over $\mathbb{Q}$. The bounds are uniform in the sense that they only depend on the rank of the corresponding Jacobian.

\begin{center}
\textbf{1. INTRODUCTION}
\end{center}
Let $F(X_0,X_1,X_2)\in \mathbb{Z}[X_0,X_1,X_2] $ be a non-singular cubic form, so that $F=0$ defines a smooth plane cubic curve $C$ in $\mathbb{P}^2.$ We want to study the asymptotic behaviour of the counting function
$$N(B)=\sharp \{P\in C(\mathbb{Q}): H(P)\leq B\},$$
with respect to the naive height function $H(P):=$max$\{|x_0|,|x_1|,|x_2|\}$ for $P=[x_0,x_1,x_2]$ with co-prime integer values of $x_0,x_1,x_2.$

\text{}\\
It is known that if the rank $r$ of the Jacobian Jac$(C)$ is positive, then we have
\begin{equation}\label{1}
N(B)\sim c_F(\text{log }B)^{r/2}
\end{equation}
as $B\rightarrow \infty.$ This result was shown by Néron. Moreover, if $r=0$ then $N(B)\leq 16$ by Mazur's theorem (see Mazur [7], Theorem 8) on torsion groups of elliptic curves. But (1) is not a uniform upper bound as the constant $c_F$ depends on $C$. The aim of this paper is to give uniform upper bounds for $N(B)$ which only depend on the rank of Jac($C$).

\text{}\\
In this direction, Heath-Brown and Testa (see [6], Corollary 1.3 ) established the \linebreak uniform bound
\begin{equation}\label{1a}
N(B)\ll (\text{log }B)^{3+r/2}
\end{equation}
by using the $p-$adic determinant method developed by the first author (see [5]). In [6], they also used a result of David [2] about the successive minima of the quadratic form given by the canonical height pairing on Jac($C$) to prove the sharper uniform bounds $N(B)\ll (\text{log }B)^{1+r/2}$ for all $r$ and $N(B)\ll (\text{log }B)^{r/2}$ if $r$ is sufficiently large.

\text{}\\
We shall in this paper give a direct proof of the bound
\begin{equation}\label{1b}
N(B)\ll (\text{log }B)^{2+r/2},
\end{equation}
based on the determinant method, which does not depend on any deep result about the canonical height pairing.

\text{}\\
To do this, we follow the approach in [6] with descent. But we replace the $p-$adic determinant method by a global determinant method developed by Salberger [9]. The main result of this paper is the following

\text{}\\
\textbf{Theorem 1.} \emph{Let $F(X_0,X_1,X_2)\in \mathbb{Z}[X_0,X_1,X_2]$ be a non-singular cubic form, so that $F=0$ defines a smooth plane cubic curve $C$. Let $r$ be the rank of Jac$(C)$. Then for any $B\geq 3$ and any positive integer $m$ we have}
$$\displaystyle N(B)\ll m^r \left(B^{\frac{2}{3m^2}}+m^2\right)\log B$$
\emph{uniformly in $C$, with an implied constant independent of} $m.$

\text{}\\
This bound improves upon the estimate
$$N(B)\ll m^{r+2}\left(B^{\frac{2}{3m^2}}\log B+\log ^2 B \right)$$
in [6] (see Theorem 1.2). Taking $m=1+[\sqrt{\log B}]$ we immediately obtain the following result.

\text{}\\
\textbf{Corollary 2.} \emph{Under the conditions above we have}
$$N(B)\ll (\text{log }B)^{2+r/2}$$
\emph{uniformly in} $C$.

\text{}\\
In the appendix we include for comparison a short account of the bounds for $N(B)$ that can be deduced from David's result.

\text{}\\
\begin{center}
\textbf{2. THE DESCENT ARGUMENT}
\end{center}
We shall in this section recall the argument in [6], where the study of $N(B)$ is reduced to a counting problem for a biprojective curve.

\text{}\\
Let $\psi : C\times C \rightarrow \text{Jac}(C)$ be the morphism to the Jacobian of $C$ defined by \linebreak $\psi(P,Q)=[P]-[Q].$ Let $m$ be a positive integer and define an equivalence relation on $C(\mathbb{Q})$ as follows: $P\sim_m Q$ if $\psi(P,Q)\in m(\text{Jac}(C)(\mathbb{Q})).$ The number of equivalence classes is at most $16m^r$ by the theorems of Mazur and Mordell-Weil. There is therefore a class $K$ such that
$$N(B)\ll m^r \sharp \{P\in K : H(P)\leq B\}.$$
If we fix a point $R$ in $K$ then for any other point $P$ in $K$, there will be a further point $Q$ in $C(\mathbb{Q})$ such that $[P]=m[Q]-(m-1)[R]$ in the divisor class group of $C$. We define the curve $X=X_R$ by
$$X_R:=\{(P,Q)\in C\times C : [P]=m[Q]-(m-1)[R]\}$$
in $\mathbb{P}^2\times \mathbb{P}^2.$ Then $N(B)\ll m^r \sharp \mathcal{K}$, where
$$\mathcal{K}:= \{(P,Q)\in X(\mathbb{Q}): H(P)\leq B\}.$$
We have thus reduced the counting problem for $C$ to a counting problem for a \linebreak biprojective curve $X$ in $\mathbb{P}^2\times \mathbb{P}^2.$ We shall also need the following lemma from [6] (see Lemma 2.1).

\text{}\\
\textbf{Lemma 3.} \emph{Let $C$ be a smooth plane cubic curve defined by a primitive form $F$ with $\|F\|\ll B^{30},$ and $R$ be a point in $C(\mathbb{Q}).$ There exists an absolute constant $A$ with the following property. Suppose that $(P,Q)$ is a point in $X_R(\mathbb{Q})$ and that $B\geq 3.$ Then if $H(P),H(R)\leq B$ we have} $H(Q)\leq B^A.$

\text{}\\
\begin{center}
\textbf{3. THE GLOBAL DETERMINANT METHOD}
\end{center}
We shall in this section apply Salberger's global determinant method in [9] to $X$ and consider congruences between integral points on $X$ modulo all primes of good reduction for $C$ and $X$. It is a refinement of the $p-$adic determinant method used in [5] and [6].

\text{}\\
We will label the points in $\mathcal{K}$ as $(P_j,Q_j)$ for $1\leq j\leq N$, say, and fix integers $a,b\geq 1$. Let $I_1$ be the vector space of all bihomogeneous forms in $(x_0,x_1,x_2;y_0,y_1,y_2)$ of bidegree $(a,b)$ with coefficients in $\mathbb{Q}$ and $I_2$ be the subspace of such forms which vanish on $X$. Since the monomials
$$x_0^{e_0}x_1^{e_1}x_2^{e_2}y_0^{f_0}y_1^{f_1}y_2^{f_2}$$
with $$e_0+e_1+e_2=a \text{  and  }f_0+f_1+f_2=b$$
form a basis for $I_1$, there is a subset of monomials $\{F_1,...,F_s\}$ whose corresponding cosets form a basis for $I_1/I_2$. As in [6] (see Lemma 3.1), if $\frac{1}{a}+\frac{m^2}{b}<3$, then \linebreak $s=3(m^2a+b).$ Thus we shall always assume that $a\geq 1$ and $b\geq m^2$ to make sure that $s=3(m^2 a +b).$ Consider the $N\times s$ matrix
$$M=\left(
  \begin{array}{cccc}
    F_1(P_1,Q_1) & F_2(P_1,Q_1) & \ldots & F_s(P_1,Q_1) \\
    F_1(P_2,Q_2) & F_2(P_2,Q_2) & \ldots & F_s(P_2,Q_2) \\
    \vdots & \vdots & \ldots & \vdots \\
    F_1(P_N,Q_N) & F_2(P_N,Q_N) & \ldots & F_s(P_N,Q_N) \\
  \end{array}
\right).$$
If we can choose $a$ and $b$ such that rank$(M)<s$, then there is a non-zero column vector $\underline{c}$ such that $M\underline{c}=\underline{0}$. This will produce a bihomogeneous form $G$, say, of bidegree $(a,b)$ such that $G(P_j,Q_j)=0$ for all $1\leq j\leq N.$ Hence all points in $\mathcal{K}$ will lie on the variety $Y\subset \mathbb{P}^2\times \mathbb{P}^2$ given by $G=0,$ while the irreducible curve $X$ does not lie on $Y$. Thus
\begin{equation}
N\leq \sharp (X\cap Y)\leq 3(m^2a+b)
\end{equation}
by the Bezout-type argument in [6] (see Lemma 5.1).

\text{}\\
In order to show that rank$(M)<s$, we may clearly suppose that $N\geq s$. We will show that each $s\times s$ minor det$(\Delta)$ of $M$ vanishes. Without loss of generality, let $\Delta$ be the $s\times s$ matrix formed by the first $s$ rows of $M$.

$$\Delta=\left(
  \begin{array}{cccc}
    F_1(P_1,Q_1) & F_2(P_1,Q_1) & \ldots & F_s(P_1,Q_1) \\
    F_1(P_2,Q_2) & F_2(P_2,Q_2) & \ldots & F_s(P_2,Q_2) \\
    \vdots & \vdots & \ldots & \vdots \\
    F_1(P_s,Q_s) & F_2(P_s,Q_s) & \ldots & F_s(P_s,Q_s) \\
  \end{array}
\right).$$
The idea is now to give an upper bound for det$(\Delta)$ which is smaller than a certain integral factor of det$(\Delta)$. To do this, we first recall a result from [5] (see Theorem 4).

\text{}\\
\textbf{Lemma 4}. \emph{For a plane cubic curve $C$ defined by a primitive integral form $F$, either} $N(B)\leq 9$ \emph{or} $\| F\| \ll B^{30}.$

\text{}\\
Thus from now on, we may and shall always suppose that $\| F\| \ll B^{30}.$ It is not difficult to see that every entry in $\Delta$ has modulus at most $B^a B^{Ab}$, where $A$ is the absolute constant in Lemma 3. Since $\Delta$ is a $s\times s$ matrix, we get that
\begin{equation}\label{2}
  \text{log}|\text{det}(\Delta)|\leq s\text{log }s+s\text{log }B^{a+Ab}.
\end{equation}
Now we find a factor of det$(\Delta)$ of the form $p^{N_p},$ where $p$ is a prime of good reduction for $C$. In order to do that, we divide $\Delta$ into blocks such that elements in each block have the same reduction modulo $p$.

\text{}\\
Let $p$ be a prime number and $Q^*$ be a point on $C(\mathbb{F}_p).$ Then we define the set
$$S(Q^*,p,\Delta )=\{(P_j,Q_j): 1\leq j \leq s, \text{ }\overline{Q_j}=Q^*\},$$
where $\overline{Q_j}$ denotes the reduction from $C(\mathbb{Q})$ to $C(\mathbb{F}_p)$. Suppose $\sharp S(Q^*,p,\Delta )=E$. We consider any $E\times E$ sub-matrix $\Delta^*$ of $\Delta$ corresponding to $S(Q^*,p,\Delta )$ and recall a result from [6] (see Lemma 4.2). Note that our set $S(Q^*,p,\Delta )$ has less elements than the set $S(Q';p,B)$ defined at the beginning of Section 3 in [6] but the proof still works.

\text{}\\
\textbf{Lemma 5.} \emph{If $p$ is a prime of good reduction for $C$. Then $p^{E(E-1)/2}$ divides} det$(\Delta^*).$

\text{}\\
From this lemma we obtain a factor of det$(\Delta)$ of the form $p^{N_p}$ by means of Laplace expansion. Moreover, we can do the same argument for all primes of good reduction for $C$ and then obtain a very large factor of det$(\Delta).$ That is the idea of the global determinant method in [9].

\text{}\\
\textbf{Lemma 6.} \emph{Let $p$ be a prime of good reduction for $C$. There exists a non-negative integer} $N_p\geq \frac{s^2}{2n_p}+ O(s)$ \emph{such that} $p^{N_p} | \text{det}(\Delta),$ \emph{where $n_p$ is the number of $\mathbb{F}_p-$points on} $C(\mathbb{F}_p)$.

\text{}\\
\emph{Proof.} Let $P$ be a point on $C(\mathbb{F}_p)$ and $s_P$ be the number of elements in $S(P,p,\Delta ).$ Then by Lemma 5, there exists an integer $N_P=s_P(s_P-1)/2$ such that $p^{N_P} | \text{det}(\Delta^*)$ for each $s_P\times s_P$ sub-matrix $\Delta^*$ of $\Delta$ corresponding to $S(P,p,\Delta )$.

\text{}\\
If we apply this to all points on $C(\mathbb{F}_p)$ and use Laplace expansion, then we get that $p^{N_p} | \text{det}(\Delta)$ for
$$N_p=\sum_P N_P=\frac{1}{2}\sum_P {s_P}^2-\frac{s}{2}\geq \frac{s^2}{2n_p}+O(s)$$
in case $C$ has good reduction at $p$. This completes the proof of Lemma 6.

\text{}\\
We now give a bound for the product of primes of bad reduction for $C$. Since \linebreak $\|F\|\ll B^{30}$, the discriminant $D_F$ of $F$ will satisfy log$|D_F|\ll$ log $B$. Thus \linebreak log $\Pi_C\ll$ log $B$, where $\Pi_C$ is the product of all primes of bad reduction for $C$. We have therefore the following bound.

\text{}\\
\textbf{Lemma 7.} \emph{Suppose that} $\|F\|\ll B^{30}.$ \emph{The product $\Pi_C$ of all primes of bad reduction for} $C$ \emph{satisfies} log $\Pi_C=O(\text{log }B).$

\text{}\\
We need one more lemma from [9] (see Lemma 1.10).

\text{}\\
\textbf{Lemma 8.} \emph{Let $\Pi>1$ be an integer and $p$ run over all prime factors of $\Pi$. Then}
$$\displaystyle \sum_{p|\Pi}\frac{\text{log }p}{p}\leq \text{log log }\Pi+2.$$\\
\emph{Proof.} We may and shall assume that $\Pi$ is a square-free. Let $l$ be a positive integer such that $l\leq \Pi$ and $v_p(n)$ be the highest integer such that $p^{v_p(n)}| n.$ We then have (see Tenenbaum [10], p. 13-14)
$$\displaystyle l\sum_{p|\Pi}\frac{\text{log }p}{p}-\sum_{p|\Pi}\text{log }p\leq \sum_{p|\Pi}v_p(l!)\text{log }p$$
$$\displaystyle\leq \sum_{p\leq\Pi}v_p(l!)\text{log }p=\text{log }l!\leq l\text{log }l,$$
$$\Rightarrow \displaystyle \sum_{p|\Pi}\frac{\text{log }p}{p}\leq \text{log }l+\frac{1}{l}\sum_{p|\Pi}\text{log }p\leq \text{log }l+(1/l)\text{log }\Pi.$$
To obtain the assertion, let $l=[\text{log }\Pi]$ for $\Pi>2.$
\newpage
\begin{center}
\textbf{4. PROOF OF THEOREM 1}
\end{center}
We now use the lemmas in Section 3 to prove that det$(\Delta)$ vanishes if $s$ is large enough. Let $\Pi_C$ be the product of all primes $p$ of bad reduction for $C$. Then

\begin{equation}\label{3}
\displaystyle \sum_{p|\Pi_C}\frac{\text{log }p}{p}\leq \text{log }\text{log }B+O(1)
\end{equation}
by Lemma 7 and Lemma 8. We apply Lemma 6 to the primes $p\leq s$ of good reduction for $C$ and write $\displaystyle {\sum_{p\leq s}}^*$ for a sum over these primes. We then obtain a positive factor $T$ of det$(\Delta)$ which is relatively prime to $\Pi_C$ such that
$$\displaystyle\text{log }T\geq \frac{s^2}{2}{\sum_{p\leq s}}^* \frac{\text{log }p}{n_p}+O(s){\sum_{p\leq s}}^* \text{log }p.$$
The last term is $O(s^2)$ since $\sum_{p\leq s}\text{log }p=O(s)$ (see [10], p. 31). Also,
$$\displaystyle \frac{\text{log }p}{n_p}\geq \frac{\text{log }p}{p}-\frac{(n_p-p)\text{log }p}{p^2}.$$
Moreover, it is well-known that if $p$ is a prime of good reduction for $C$, then \linebreak $n_p=p+O(\sqrt{p}).$ Thus we conclude that
$$\displaystyle \frac{\text{log }p}{n_p}\geq \frac{\text{log }p}{p}+O\left(\frac{\text{log }p}{p^{3/2}}\right)$$
for all primes $p$ of good reduction for $C.$ Therefore,
$${\sum_{p\leq s}}^* \frac{\text{log }p}{n_p}\geq {\sum_{p\leq s}}^* \frac{\text{log }p}{p}+O(1)$$
and then $$\text{log }T\geq \frac{s^2}{2}{\sum_{p\leq s}}^* \frac{\text{log }p}{p}+O(s^2)$$
But by (6),
$$\displaystyle\sum_{p\leq s} \frac{\text{log }p}{p}-{\sum_{p\leq s}}^* \frac{\text{log }p}{p}\leq \text{log }\text{log }B+O(1)$$
and $\displaystyle \sum_{p\leq s} \frac{\text{log }p}{p}=\text{log }s+O(1)$ (see [10], p. 14). Hence,
\begin{equation}\label{4}
\displaystyle \text{log }T\geq \frac{s^2}{2}\text{log}\left(\frac{s}{\text{log }B}\right)+O(s^2).
\end{equation}
Thus from (5) and (7) we obtain
$$ \displaystyle \text{log}\left(\frac{|\text{det}(\Delta)|}{T}\right)\leq s\text{log }s+s\text{log }B^{a+Ab}-\frac{s^2}{2}\text{log}\left(\frac{s}{\text{log }B}\right)+O(s^2)$$
$$=\displaystyle \frac{s^2}{2}\left(\log B^{\frac{2(a+Ab)}{s}}-\text{log}\left(\frac{s}{\text{log }B}\right)\right)+O(s^2).$$
There is therefore an absolute constant $u\geq 1$ such that
$$\displaystyle \text{log}\left(\frac{|\text{det}(\Delta)|}{T}\right)\leq \frac{s^2}{2}\left(\log B^{\frac{2(a+Ab)}{s}}-\text{log}\left(\frac{s}{u\text{log }B}\right)\right).$$
If
\begin{equation}
\displaystyle s>u B^{\frac{2(a+Ab)}{s}}\text{log }B
\end{equation}
we have in particular that $\log \left(\frac{|\det(\Delta)|}{T}\right)<0$ and hence det$(\Delta)=0$ as $\frac{|\det(\Delta)|}{T}\in \mathbb{Z}_{\geq 0}.$

\text{}\\
Remember that $s=3(m^2 a +b)$ if $a\geq 1$ and $b\geq m^2.$ We now choose $b=m^2$ and
$$a=1+\left[\frac{uB^{\frac{2}{3m^2}}\log B}{m^2}+A\log B\right].$$
Then
$$u B^{\frac{2(a+Ab)}{s}}\text{log }B=u B^{\frac{2(a+Am^2)}{3m^2(a+1)}}\text{log }B $$
$$< u B^{\frac{2}{3m^2}}B^{\frac{2A}{3a}}\log B<s.$$
Thus (8) holds and hence det$(\Delta)=0$. Then rank$(M)<s$ such that there is a \linebreak bihomogeneous form in $\mathbb{Q}[x_0,x_1,x_2,y_0,y_1,y_2]$ which vanishes at all $(P_j,Q_j)\in X(Q), \linebreak \text{ }1\leq j\leq N,$ with $H(P_j)\leq B$ but not everywhere on $X$. Hence (see (4))
$$N\leq 3(m^2 a+b)\ll \left(B^{\frac{2}{3m^2}}+m^2\right)\log B$$
$$\Rightarrow N(B)\ll m^r \left(B^{\frac{2}{3m^2}}+m^2\right)\log B.$$
This completes the proof of Theorem 1.

\text{}\\
\textbf{Acknowledgement}\\
I wish to thank my supervisor Per Salberger for introducing me to the problem and giving me many useful suggestions.

\text{}\\
\begin{center}
\textbf{APPENDIX}
\end{center}
In this appendix we record the following more precise version of a result in [6].

\text{}\\
\textbf{Theorem 9.} \emph{Let $C$ be any smooth plane cubic curve and $r$ be the rank of Jac($C$). Let $m_l=\frac{l^2-4l-4}{8l^2+8l}$ for $l\geq 1$. Then}
$$N(B)\ll \left\{
            \begin{array}{ll}
              (\log B)^{-(m_1+...+m_r)+r/2}, & \hbox{if $1\leq r <16$;} \\
              (\log B)^{r/2}, & \hbox{if $r\geq 16$.}
            \end{array}
          \right.$$
\emph{with an absolute implied constant. In particular,} $N(B)\ll (\log B)^{1+r/2}$ \emph{for all} $r$.

\text{}\\
\emph{Proof.} The proof is just a careful re-examination of the argument of Heath-Brown and Testa [6]. This argument is based on a result of David [2] about successive minima for the quadratic form $Q$ corresponding to the canonical height on Jac($C$). As in [6] (see (11)),
\begin{equation}
N(B)\ll \displaystyle \prod_{j\leq r}\text{max}\left\{1,4\frac{\sqrt{c\text{log }B}}{M_j}\right\},
\end{equation}
where $c$ is an absolute constant and $M_j,$ $j=1,...,r$ are successive minima of $\sqrt{Q}$.

\text{}\\
We now recall Corollary 1.6 from [2], which shows that if $D$ is the discriminant of Jac($C$) then for all $l\leq r$, $M_l\gg (\text{log}|D|)^{m_l},$ where $m_l=\frac{l^2-4l-4}{8l^2+8l}.$ Note that David's result refers to the successive minima for $Q$, while we have given the corresponding results for $\sqrt{Q}.$

\text{}\\
In Lemma 4 we saw that $\|F\|\ll B^{30}$ if $N(B)>9.$ There is, therefore, in that case an absolute constant $k$ such that
$$\max \left\{1,4\frac{\sqrt{c\log B}}{M_j}\right\}\leq k (\log B)^{1/2}(\log |D|)^{-m_j}$$
for $j=1,...,r$ since $|m_j|<1/2$ and $\log |D|\ll \log B.$ Hence, if $N(B)>9$, then from (9) we obtain
\begin{equation}
N(B)\ll k^r (\log B)^{r/2} (\log |D|)^{-(m_1+...+m_r)}.
\end{equation}
If $1\leq r<16,$ then $-(m_1+...+m_r)>0$ and the assertion holds. If $r\geq 16,$ let $D_0=\exp(k^{1/m_{16}}).$ Then $k(\log |D|)^{-m_j}\leq 1$ for $j>16$ and $|D|\geq D_0.$ Hence
$$N(B)\ll (\log B)^{r/2} (\log |D|)^{-(m_1+...+m_{16})}\ll (\log B)^{r/2}$$
as $-(m_1+...+m_{16})<0.$ When $|D|\leq D_0$ the rank $r$ is bounded and we get the same assertion by (10).

\text{}\\
So in any case, $N(B)\ll (\text{log }B)^{r/2},$ if $r\geq 16.$ It should thereby be noted that Elkies (see [3]) has shown that there exist elliptic curves of rank $r \geq 28$.

\end{document}